\documentclass[12pt]{amsart}
\usepackage{amsmath,amsthm,amsfonts,amssymb,eucal}

\renewcommand {\a}{ \alpha }

\newcommand{\vare}{\varepsilon}

\newcommand{\varf}{\varphi}

\newcommand{\D}{\Delta}

\renewcommand{\l}{\lambda}
\renewcommand{\L}{\Lambda}

\newcommand{\T}{\Theta}
\newcommand{\p}{\partial}

\newcommand{\R}{ \mathbb R}

\newcommand{\N}{ \mathbb N}

\newcommand {\bb}{\mathbf b}

\newcommand {\BB}{\mathbf B}

\newcommand {\BG}{\mathbf G}

\newcommand {\BK}{\mathbf K}
\newcommand {\BL}{\mathbf L}

\newcommand {\BT}{\mathbf T}

\newcommand{\lu}{\langle}
\newcommand{\ru}{\rangle}

\newcommand{\CV}{\mathcal V}

\newcommand{\CL}{\mathcal L}

\newcommand{\CP}{\mathcal P}

\newcommand{\cm}{\mathcal M}

\newcommand{\CE}{\mathcal E}

\newcommand{\tC}{\textsf C}

\newcommand{\tH}{\textsf H}

\newcommand{\tL}{\textsf L}

\newcommand{\tS}{\textsf S}

\newcommand{\tW}{\textsf W}

\DeclareMathOperator {\diam} {{diam}}

\DeclareMathOperator{\res}{\restriction}

\newtheorem{thm}{Theorem}[section]
\newtheorem{cor}[thm]{Corollary}
\newtheorem{lem}[thm]{Lemma}

\theoremstyle{definition}

\theoremstyle{remark}

\numberwithin{equation}{section}

\newcommand{\thmref}[1]{Theorem~\ref{#1}}

 \DeclareMathOperator {\Step}{Step}
 
\DeclareMathOperator {\comb}{comb} \DeclareMathOperator {\Int}{Int}
\newcommand {\wt}{\widetilde} 
 
\begin{document} 
\hoffset -4pc 
\title[Laplacian on metric graphs] 
{On the eigenvalue estimates for the weighted Laplacian on metric graphs}
\author[M. Solomyak] 
{Michael Solomyak} 
\address 
{Department of Mathematics 
\\ Weizmann Institute\\ Rehovot\\ Israel}
\email{solom@wisdom.weizmann.ac.il} 
\thanks
{The work was supported in part by the Minerva center for non-linear 
physics and by the Israel Science Foundation.} 
\date{28 December 2001}
\dedicatory
{Dedicated to Olga
Aleksandrovna Ladyzhenskaya\\ to whom I am indebted for my formation as a
mathematician}

\begin{abstract} It is shown that the eigenvalues of the equation $-\l\D
u=Vu$ on a graph $\BG$ of final total length $|\BG|$, with non-negative
$V\in \tL^1(\BG)$ and under appropriate boundary conditions, satisfy the
inequality $n^2\l_n\le |\BG|\int_{\BG} Vdx$, independently of geometry of
a given graph. Applications and generalizations of this result are also
discussed. 
\end{abstract} 
\maketitle 

\section{Introduction}
 
Differential operator on a metric graph $\BG$ is a family of
differential expressions on its edges, complemented by appropriate
matching conditions at the most of vertices and by some boundary
conditions at the remaining vertices. In particular, for the Laplacian
$\D$ the differential expression is $\D u=u''$ and the matching conditions
at vertices are the Kirchhoff conditions coming from the theory of
electric networks. Our main goal in this paper is to investigate the
behaviour of eigenvalues for the problem 
\begin{equation}\label{1:1}
-\l\D u=Vu;\qquad u(x_0)=0,\ u'(v)=0\ {\text{if}}\
v\in\p\BG\setminus\{x_0\} 
\end{equation} 
on graphs of finite total length
$|\BG|$. In \eqref{1:1} $V$ is a given real-valued, measurable weight
function on $\BG$ and $x_0\in\BG$ is a selected vertex. Note that the set
$\p\BG\setminus\{x_0\}$ can be empty. 
 
It is possible also to consider the similar problem with zero boundary
conditions at several vertices, $u(x_1)=\ldots=u( x_r)=0$, however this
gives nothing new: by the variational principle, any estimate for the
problem \eqref{1:1} implies the same estimate for the $r$-point problem. 
The accurate setting of the problem uses quadratic forms, see Section 3
for details.  The way to insert the spectral parameter in \eqref{1:1} is
motivated by technical reasons. 
 
\vskip0.2cm
 
Here is our main result on the problem \eqref{1:1}. Let $\pm\l_n^{\pm}$
denote the positive and the negative eigenvalues of this problem. 
 
\begin{thm}\label{1:eig} 
Let $\BG$ be a connected graph of finite total
length, $x_0\in\BG$ be its arbitrary vertex, and let $V=\overline
V\in\tL^1(\BG)$. Then the eigenvalues of the problem \eqref{1:1} satisfy
the inequality 
\begin{equation}\label{1:ei} 
\l_n^{\pm}\le
\frac{|\BG|\int_{\BG} V_{\pm}dx}{n^2},\qquad \forall n\in\N 
\end{equation}
where $2V_{\pm}=|V|\pm V$. Along with the estimate \eqref{1:ei}, the
Weyl-type asymptotics holds: 
\begin{equation}\label{1:3}
n\sqrt{\l_n^{\pm}}\to\pi^{-1}\int_{\BG}\sqrt{V_{\pm}(x)}dx,\qquad
n\to\infty. 
\end{equation} 
\end{thm}
 
\vskip0.2cm 
The simplest example of metric graph is a single segment
$[0,L]\subset\R$. The basic results for this case were obtained by M.Sh. 
Birman and the author \cite{BS0} as far back as in 1971, as a particular
case of the general result for many dimensions; see also an exposition in
\cite{BS2}. Namely, it was shown that eigenvalues of the equation $-\l
u''=Vu$ on the interval $(0,L)$, under the Dirichlet boundary conditions
at its ends, satisfy the asymptotics \eqref{1:3} and the inequality
$n^2\l_n^{\pm}\le C|\BG|\int_{\BG} V_{\pm}dx$, with some absolute constant
$C$. The problem of the best possible
constant in the inequality was not discussed in \cite{BS0} and in
\cite{BS2}, though the sharp estimate with $C=1/4$ could be easily derived from
Theorem 2.2 of the paper \cite{BS1}. For the zero boundary condition
only at one point, as in \eqref{1:1}, the constant $C=1$ 
in \eqref{1:ei} is unimprovable even for the case of segment.

The most important feature of the estimate \eqref{1:ei} for $\BG=[0,L]$
is its uniformity
with respect to $V\in \tL^1(0,L)$. This proved quite useful for various
applications. \thmref{1:eig} shows that this estimate extends to arbitrary
graphs of finite total length, with the same constant as for the single
interval. The very possibility of such extension might look problematic, 
since the
eigenvalues depend on the combinatorial structure of the graph, which can
be quite diverse.

One can give a qualitative explanation of this effect. If $V\ge 0$, the
eigenvalues of the problem \eqref{1:1} can be defined in terms of
approximative characteristics of the unite ball of the Sobolev space
$\tH^1(\BG,x_0)$ as a compact in the weighted space $\tL^2(\BG,V)$;
indication of the point $x_0$ in the first notation reflects the boundary
condition $u(x_0)=0$. The points
of a graph $\BG$ lie ``closer to each other'' than the points of the
segment $[0,|\BG|]$.  Therefore, one should expect that the dispersion of
the values of a given function $u\in \tH^1 (\BG,x_0)$ is smaller than that of
a function $v\in \tH^1 (0,|\BG|)$ having the same $\tL^2$-norm of the
derivative. As a consequence, it is easier to approximate the class 
$\tH^1(\BG,x_0)$ than $\tH^1 (0,|\BG|)$. Eventually this leads to the 
estimate \eqref{1:ei}.

This, a little bit naive argument gives no clue to the proof. 
Our proof is based
upon a result (\thmref{2:par}) of a rather combinatorial nature. This
theorem can be considered as a far going generalization of Theorem 2.2
from \cite{BS1}.

Our results should be compared with the ones of the recent paper \cite{EHL} 
by W.D. Evans, D.J. Harris, and J. Lang. Its subject is the behaviour
of approximation numbers of the Hardy-type integral operators in 
the spaces $\tL^p$
on trees. If $p=2$, the approximation numbers coincide with the singular
numbers, and this is the link between our corresponding results. 
Some of our results, including the asymptotics \eqref{1:3},
could be derived from this paper. The techniques of \cite{EHL},
see especially its Section 3, is also based upon a 
certain combinatorial construction. It is substantially different from our one.
The applications discussed in our paper are not touched upon by the authors
of \cite{EHL}.

\vskip0.2cm

Let us describe the structure of the paper. In the next Section 2 we give
the necessary information about graphs and trees and prove \thmref{2:par}
which is our main technical tool. 
In Section 3 we define the Sobolev
space $\tH^1(\BG,x_0)$ and  give
the variational formulation of the problem \eqref{1:1}. We also state
two auxiliary results, Theorems \ref{3:w} and \ref{3:wid}, on approximation 
on graphs. They are useful for
the proof of \thmref{1:eig}.
The proofs of all three theorems are given in Section 4.
In Section 5 we discuss the result and give its application to the estimates 
of singular numbers for integral operators on graphs. In the last Section 6
we present a higher order analog of \thmref{1:eig}.
\vskip0.2cm

The Hilbert space structure is
unnecessary for applications of \thmref{2:par}. It can be also applied to
piecewise-polynomial approximation of Sobolev spaces $\tW^{l,p}$ on graphs
and trees. This leads to estimates of approximation numbers of embedding
of $\tW^{l,p}$ in the weighted spaces $\tL^p(\BG,V)$ on graphs of finite
total length. The results can be also applied to graphs and trees of 
infinite total length, in the spirit of the papers
\cite{NS} and \cite{EHL}.
This material will be presented elsewhere.

\section{Graphs and trees. Partitions of a tree}
 
Let $\BG$ be a graph with the set of vertices $\CV=\CV(\BG)$ and the set
of edges $\CE=\CE(\BG)$. Each edge $e$ of a metric graph is viewed as a
non-degenerate line segment of finite length $|e|$. The quantity
$|\BG|=\sum\limits_{e\in\CE(\BG)}|e|$ 
is called the
{\sl total length} of the graph $\BG$. In this paper we always assume
$|\BG|<\infty$. Even though, the number of edges of $\BG$ can be infinite.
The distance
$\rho(x,y)$ between any two points $x,y\in\BG$, and thus the metric
topology on $\BG$, is introduced in a natural way. The natural measure
$dx$ on $\BG$ is induced by the Lebesgue measure on its edges. For a
measurable set $E\subset\BG$, its measure is denoted by $|E|$. This is
consistent with the above notations $|e|,\ |\BG|$. 
 
We always consider connected graphs, unless otherwise stipulated.  We do
not exclude graphs with multiple joins. In order to avoid unnecessary 
complications, we exclude graphs with loops. Recall that a loop is 
an edge whose endpoints coincide with each other.  

For vertices $v,w$ the
notation $w\sim v$ means that there exists an edge $e\in\CE$ whose endpoints
are $v$ and $w$. Connectedness of the graph means that for any two
different vertices $v,w\in\CV$ there exists a finite sequence
$\{v_k\}_{0\le k\le m}$ of vertices, such that $v_0=v,\ v_m=w$ and
$v_k\sim v_{k-1}$ for each $k=1,\ldots,m$. The {\sl combinatorial
distance} $\rho_{\comb}(v,w)$ is defined as the minimal possible $m$ in
this construction. We define $\rho_{\comb}(v,v)=0$ for any
$v\in\CV$. The {\sl degree} $d(v)$ of a vertex $v$ is defined as the total
number of edges incident to $v$. The vertices $v$ with $d(v)=1$ constitute
the boundary $\p\BG$ of the graph $\BG$. We suppose that $d(v)<\infty$ for
each $v\in\CV$. 
It is often convenient to treat an arbitrary point $x\in\BG$ as a vertex.
We set $d(x)=2$ for any $x\notin\CV(\BG)$ and write $v\sim x$ if $v\in\CV(\BG)$
is one of the endpoints of the vertex containing $x$. This remark concerns
also the point $x_0$ appearing in \eqref{1:1} which actually can be
an arbitrary point in $\BG$.

Given a subgraph $G\subset\BG$, we denote by $d_G(v)$ the degree of a
vertex $v$ with respect to $G$. 
 
A metric graph $\BG$ is compact if and only if $\#\CE(\BG)<\infty$. 
Any metric graph can be represented as the union of an expanding family
of its compact subgraphs, 
\begin{equation}\label{2:su}
\BG=\bigcup_{m=1}^\infty \BG^{(m)},\qquad \BG^{(1)}\subset
\BG^{(2)}\subset\ldots;\qquad \#\CE(\BG^{(m)})<\infty,\ \forall m\in\N. 
\end{equation} 
Indeed,
choose an arbitrary vertex $x_0\in\BG$ and define the subgraph $\BG^{(m)}$
as follows. A vertex $v\in\CV(\BG)$ belongs to $\CV(\BG^{(m)})$ if and only if
$\rho_{\comb}(x_0,v)\le m$, and an edge $e\in\CE(\BG)$ belongs to
$\CE(\BG^{(m)})$ if and only if its both ends lie in $\CV(\BG^{(m)})$. By the
construction, the graph $\BG^{(m)}$ is connected and compact.
It is evident that $\{\BG^{(m)}\}$
is an expanding family which covers the whole of $\BG$.  
\vskip0.2cm
 
let $\BG$ be a compact graph and $G,G_1,\ldots, G_n$ be its (connected) 
subgraphs.  We say that the subgraphs $G_1,\ldots, G_n$ constitute a {\sl
partition}, or a {\sl splitting} of $G$ if $G_1\cup\ldots\cup G_n=G$ and
$|G_i\cap G_j|=0$ for any $i,j\in\{1,\ldots,n\},\ i\neq j$.  If $G,G_1$ are
subgraphs of $\BG$ and $G_1\subset G$, then connected subgraphs
$G_2,\ldots, G_n$ can be always found which together with $G_1$ constitute
a partition of $G$ (the property of {\sl complementability}). 
 
A pair $\{G,x\}$ where $G\subset\BG$ is a subgraph and $x\in G$ is a 
selected
point, is called a {\sl punctured subgraph}. If $\{G_j,x_j\},\
j=1,\ldots,n$ are punctured subgraphs such that $G=G_1\cup\ldots\cup G_n$
is a partition of a subgraph $G$, then we say that $G$ is split into the
union of punctured subgraphs. 
 
\vskip0.2cm
 
Our auxiliary result, \thmref{2:par}, which serves us as a basis for the
proof of Theorem \ref{1:eig}, concerns trees rather than arbitrary graphs. 
Recall that tree is a connected graph without cycles, loops and 
multiple joins. So, let $\BG=\BT$ be a tree.
For any two points $x,y\in\BT$ there exists a unique simple
polygonal path  in $\BT$ connecting $x$ with $y$, we denote it by 
$\lu x, y\ru$.

What was said above about partitions of graphs, applies to trees. 
Note that if $T=T_1\cup T_2$ is a partition, then the intersection
$T_1\cap T_2$ consists of exactly one point. 
\vskip0.2cm
 
Let $\Phi$ be a non-negative function defined on the set of all
subtrees $T$ of a given tree $\BT$. We call it {\sl continuous} if
$\Phi(T)\to\Phi(T_0)$ as soon as $|T\triangle T_0|\to 0$. Recall that
$T\triangle T_0=(T\setminus T_0)\cup( T_0\setminus T) $ is the symmetric
difference of the sets $T,T_0$. 
 
We call a continuous function $\Phi$ {\sl superadditive} and write
$\Phi\in\tS(\BT)$ if for any subtree $T\subset\BT$ and any partition
$T=T_1\cup\ldots\cup T_n$ we have 
\begin{equation}\label{2:super}
\Phi(T_1)+\ldots+\Phi(T_n)\le\Phi(T).  
\end{equation} 
Due to the
complementability, \eqref{2:super} implies that any superadditive function
is monotone:  
\begin{equation*}
T_1\subset T\Longrightarrow
\Phi(T_1)\le\Phi(T). 
\end{equation*}
 
Any finite Borel measure $\mu$ on $\BT$ without atoms (i.e. points of
positive measure) generates the continuous superadditive function
$\Phi(T)=\mu(T)$.  A more general example follows from H\"older's
inequality: 
\begin{equation}\label{2:gen} 
\Phi_1,\Phi_2\in\tS(\BT),\
\a_1,\a_2>0,\ \a_1+\a_2=1\ \Longrightarrow\
\Phi_1^{\a_1}\Phi_2^{\a_2}\in\tS(\BT). 
\end{equation}

Let $\{T,x\}$ be a punctured subtree of $\BT$. The tree $T$ splits in a
unique way into the union of subtrees $\T_j\subset T$, $j=1,\ldots,
d_T(x)$, rooted at $x$ and such that $d_{\T_j}(x)=1$ for each $j$. We call
it {\sl the canonical partition} of the punctured subtree $\{T,x\}$. Given
a function $\Phi\in\tS(\BT)$, we define the function $\wt\Phi(T,x)$ of
punctured subtrees, 
\begin{equation}\label{2:punc} 
\wt\Phi(T,x)=\max_{1\le j\le d_T(x)}\Phi(\T_j). 
\end{equation} 
Evidently,
$\wt\Phi(T,x)\le\Phi(T)$. 
 
Let in particular $T=\BT$. Each subtree $\T_j$ appearing in the canonical
partition of $\{\BT,x\}$ is fully determined by indication of its initial
edge $\lu x, v\ru$, $v\sim x$ and we denote this subtree by $\T_{\lu
x,v\ru}$. For $T=\BT$ the definition \eqref{2:punc} takes the form
\begin{equation*}
\wt\Phi(\BT,x)=\max_{v\sim x}\Phi(\T_{\lu x,v\ru}). 
\end{equation*}
 
\bigskip
 
Our eigenvalue estimates will be derived from the following result on
superadditive functions of subtrees.

\begin{thm}\label{2:par} 
Let $\BT$ be a compact metric tree and
$\Phi\in\tS(\BT)$. Then for any $n\in\N$ the tree $\BT$ can be split into
the union of punctured subtrees $\{T_j,x_j\}$, $j=1,\ldots,k$ in such a
way that $k\le n$ and
 
\begin{equation}\label{1:main}
\max\limits_{j=1,\ldots,k}\wt\Phi(T_j,x_j)\le (n+1)^{-1}\Phi(\BT). 
\end{equation} 
\end{thm}

\vskip0.2cm
 
First, we prove a lemma.  
\begin{lem}\label{1:lem} 
Let $\BT$ be a compact
metric tree and $\Phi\in\tS(\BT)$. Then for any $\vare$,
$0<\vare<\Phi(\BT)$, there exists a partition $\BT=T\cup T'$, such that
\begin{equation}\label{1:03} 
\Phi(T')\le \Phi(\BT)-\vare 
\end{equation}
and for the only point $x\in T\cap T'$ the inequality holds: 
\begin{equation}\label{1:04} 
\wt\Phi(T,x)\le\vare. 
\end{equation}
\end{lem} 
\begin{proof} 
Without loss of generality, we can assume
$\Phi(\BT)=1$. Take any vertex $v_0\in\p\BT$, then
$\wt\Phi(\BT,v_0)=\Phi(\BT)=1$. There is a unique vertex $v_1\sim v_0$.
Now we choose the vertices $v_2\sim v_1,\ldots, v_{k+1}\sim v_k,\ldots$ as
follows. If $v_k$ is already chosen, we define $v_{k+1}$ as the vertex
different from $v_{k-1}$ and such that
 \begin{equation}\label{1:path} 
\Phi(\T_{\lu v_k,v_{k+1}\ru})= \max_{w\sim
v_k,w\neq v_{k-1}}\Phi(\T_{\lu v_k,w\ru}) =\wt\Phi(\T_{\lu
v_k,v_{k+1}\ru},v_k). 
\end{equation} 
If there are several vertices $w\sim
v_k$ at which the maximum in the middle term of \eqref{1:path} is
attained, then any of them can be chosen as $v_{k+1}$. The described
procedure is always finite, it terminates when we arrive at a vertex
$v_m\in\p\BT$. On the path $\CP= \lu v_0,v_m\ru$ we
introduce the natural ordering, i.e. $y\succeq x$ means that $x\in\langle
v_0,y\rangle$. We write $y\succ x$ if $y\succeq x$ and $y\neq x$. 
 
Let $x\in \CP$ be not a vertex of $\BT$, then $v_{k-1}\prec x\prec
v_k$ for some $k=1,\ldots,m$. Denote 
\begin{equation*} 
T_x^+=\T_{\lu x,v_k\ru},\qquad T_x^-=\T_{\lu x,v_{k-1}\ru}. 
\end{equation*}
We also define the subtrees $T_x^\pm$ for $x=v_0,\ldots, v_m$. Namely,
\begin{gather*} 
T_{v_k}^-=T_{\langle v_k,v_{k-1}\rangle},\qquad k=1,\ldots,
m;\\ T_{v_k}^+=\bigcap_{v_{k-1}\prec x\prec v_k}T_x^+= \bigcup_{v\sim v_k,
v\neq v_{k-1}}T_{\langle v_k,v\rangle}, \qquad k=0,\ldots, m-1.
\end{gather*} 
Finally, $T_{v_0}^-=\{v_0\},\ T_{v_m}^+=\{v_m\}$ are
degenerate subtrees. Clearly, for any $x\in\CP$
$\BT=T_x^+\cup T_x^-$ 
is a partition of the punctured tree $\{\BT,x\}$, and $T_x^+\cap
T_x^-=\{x\}$. 
\vskip0.2cm
 
The function $F(x)=\Phi(T_x^+)$ is well defined on $\CP$ and is continuous
everywhere except possibly for the vertices $v_0,\ldots, v_m$. By the
construction, 
\begin{equation*} 
\wt\Phi(T^+_x,x)=F(x),\qquad x\neq
v_0,\ldots, v_m, 
\end{equation*} and 
\begin{gather*} \lim_{x\prec v_k,x\to
v_k}F(x)=F(v_k);\\ \lim_{x\succ v_k,x\to v_k}F(x)=\Phi(T_{\langle
v_k,v_{k+1}\rangle}) =\wt\Phi(T_{v_k}^+,v_k)\le F(v_k). 
\end{gather*} 
Here the second equality follows from \eqref{1:path}. It is clear that $F(x)$
is non-increasing along the path $\CP$. Besides, 
\begin{equation*}
0=F(v_m)<\vare<F(v_0)=1. 
\end{equation*} 
Therefore, there exists a point $x\in \CP$ such that
\begin{equation*} 
\wt{\Phi}(T_x,x)\le\vare \le F(x).
\end{equation*}
We take $T=T_x^+$ and $T'=T_x^-$. The inequality \eqref{1:04} is satisfied
and \eqref{1:03} is implied by superadditivity: 
\begin{equation*}
\Phi(T')\le 1-\Phi(T)=1-F(x)\le 1-\vare. 
  \end{equation*} 
\end{proof}
 
\noindent {\it Proof of \thmref{2:par}.} 
1. Let $n=1$. Then we apply the
result of Lemma \ref{1:lem} with $\vare=\Phi(\BT)/2$. Let $\BT=T\cup T'$
be the corresponding partition, then
$\wt\Phi(T',x)\le\Phi(T')\le\Phi(\BT)/2$.  Consider the canonical
partition of the punctured tree $\{\BT,x\}$. Each subtree of this
partition is contained either in $T$ or in $T'$, therefore
\begin{equation*}
\wt\Phi(\BT,x)\le\max\bigl(\wt\Phi(T,x),
\wt\Phi(T',x)\bigr)\le\Phi(\BT)/2. 
\end{equation*} 
Thus, \eqref{1:main} with $k=n=1$ is
satisfied if we take $T_1=\BT$ and $x_1=x$. 
 
\vskip0.2cm 2. We proceed by induction. Suppose that the result is already
proved for $n=n_0-1$.  Let $\BT=T\cup T'$ be the partition constructed
according to Lemma \ref{1:lem} for $\vare=(n_0+1)^{-1}\Phi(\BT)$. Then
\begin{equation*} 
\Phi(T')\le n_0(n_0+1)^{-1}\Phi(\BT). 
\end{equation*} 
By the inductive hypothesis, there exists a splitting of
$T'$ into the union of the family of punctured subtrees $\{T_j,x_j\}$,
$j=1,\ldots,k$ such that $ k\le n_0-1 $ and for each $j$ 
\begin{equation*}
\widetilde\Phi(T_j,x_j)\le n_0^{-1}\Phi(T')\le (n_0+1)^{-1}\Phi(\BT). 
\end{equation*} 
Adding to this family the punctured subtree
$\{T_{k+1},x_{k+1}\}=\{T,x\}$, we obtain the desired partition of $\BT$ 
for $n=n_0$.\qed
 
\section{Variational setting of the problem. Reduction to the case of
trees}
 
\subsection{ Sobolev spaces on a graph.} Below $\|\cdot\|_p,\ 1\le
p\le\infty$ stands for the norm in the space $\tL^p(\BG)$ and $\tL_+(\BG)$
stands for the cone of all non-negative elements in $\tL^1(\BG)$. 
 
We say that a function $u$ on $\BG$ belongs to the Sobolev space
$\tL^{1,2}(\BG)$ if $u$ is continuous on $\BG$, the restriction of $u$
to each edge $e$ lies in $\tH^1(e)$, and $u'\in \tL^2(\BG)$. The
functional $\|u'\|_2$ defines on $\tL^{1,2}(\BG)$ a semi-norm which vanishes
on the one-dimensional subspace of constant functions. 
 
Let $\BG$ be a graph of finite total length and let $\xi,x$ be its two
arbitrary points. Choose a simple polygonal path $\CL$ in $\BG$ connecting
$\xi$ with $x$, let its length be $t_0$. Parametrizing $\CL$ by the path
length, we can regard the restriction $u\res\CL$ as a function on the line
segment $[0,t_0]$. It follows from the equality
$u(x)-u(\xi)=\int_0^{t_0}u'(t)dt$
that
\begin{equation}\label{3:m} 
|u(x)-u(\xi)|^2\le
t_0\int_0^{t_0}|u'(t)|^2dt\le |\BG|\int_{\BG}|u'(x)|^2dx. 
\end{equation}
This shows that any function $u\in\tL^{1,2}(\BG)$ lies in the H\"older class
of order $1/2$. 
\vskip0.2cm

A {\sl step function} $v$ on $\BG$ is a function which takes only a finite
number of different values, each on a connected subset of $\BG$. 
We denote by $\Step(\BG)$
the linear space (non-closed linear subspace of $\tL^\infty(\BG)$) of all
step functions on $\BG$. 
 
The following result on the approximation of functions $u\in\tL^{1,2}(\BG)$ by
step-functions will be used when proving \thmref{1:eig}. Let punctured
subgraphs $\{G_j,x_j\},\ j=1,\ldots,k$ constitute a partition of the graph
$\BG$. We associate with this partition the linear operator
\begin{equation}\label{3:op} 
P:u\mapsto v=\sum_{j=1}^k u(x_j)\chi_j
\end{equation} 
where $\chi_j$ stands for the characteristic function of
the set $G_j$.  It is clear that the operator $P$ acts from $\tL^{1,2}(\BG)$
into $\Step(\BG)$ and its rank is less or equal to $k$. 
 
\begin{thm}\label{3:w} 
Let $\BG$ be a compact graph and $V\in\tL_+(\BG)$. 
Then for any $n\in\N$ there exists a partition of $\BG$ into punctured
subgraphs $\{G_j,x_j\},\ j=1,\ldots,k$, such that $k\le n$ and for the
corresponding operator $P$ given by \eqref{3:op} we have
\begin{equation}\label{3:2} 
\int_{\BG}|u-P u|^2Vdx\le
\frac{|\BG|\int_{\BG}Vdx}{(n+1)^2} \|u'\|^2_2,\qquad \forall
u\in\tL^{1,2}(\BG). 
\end{equation} 
\end{thm}
 
The assumption that the graph $\BG$ is compact, is important for the proof. 
To exhaust the general case, we need one more statement. 
 
A graph $\BG$ of finite total length is not necessarily a compact metric
space. Let $\overline{\BG}$ be its compactification. Any function
$u\in\tL^{1,2}(\BG)$ is uniformly continuous on $\BG$ and hence admits the
unique continuous extension to $\overline{\BG}$. We keep the same
symbol $u$ for the extended function. 
 
\begin{thm}\label{3:wid} 
Let $\BG$ be a graph of finite total length and
$V\in\tL_+(\BG)$. Then for any $n\in\N$ there exist points $\overline
x_1,\ldots,\overline x_k\in \overline{\BG}$ such that $k\le n$ and the
inequality
 
\begin{equation}\label{3:wi} 
\int_{\BG}|u^2|Vdx\le
\frac{|\BG|\int_{\BG}Vdx}{(n+1)^2}\|u'\|^2_2 
\end{equation} 
holds for any
function $u\in\tL^{1,2}(\BG)$ satisfying the conditions
$u(\overline x_1)=\ldots=u(\overline x_k)=0$.  
\end{thm}
 
Proofs of Theorems \ref{3:w} and \ref{3:wid} are given in the next
section, before proving \thmref{1:eig}. 
 
\subsection{Space $\tH^1(\BG,x_0)$ and operators $\BB_V$.} 
Let a point $x_0\in\BG$ be given. It is convenient to assume that
$x_0$ is a vertex. Consider the Hilbert space 
\begin{equation*}
\tH^1(\BG,x_0)=\bigl\{u\in\tL^{1,2}(\BG):u(x_0)=0\bigr\},
\end{equation*}
equipped with the scalar product
\begin{equation*}
(u,v)_{\tH^1(\BG,x_0)}=(u',v')_{\tL_2(\BG)}.
\end{equation*}
Inequality \eqref{3:m} (with $\xi=x_0$)
shows that this scalar product is non-degenerate.
\vskip0.2cm
Let $V$ be a function  from the space $\tL^1(\BG)$. 
Consider the quadratic form 
\begin{equation}\label{3:s}
\bb_V[u]=\int_{\BG}|u|^2Vdx. 
\end{equation} 
It follows from the inequality
\eqref{3:m} (again, with $\xi=x_0$)  that  $\bb_V$ is bounded
in the space $\tH^1(\BG,x_0)$, i.e. 
\begin{equation}\label{3:r}
|\bb_V[u]|\le |\BG|\|u'\|^2_2\int_{\BG}|V(x)|dx,\qquad\forall
u\in\tH^1(\BG,x_0). 
\end{equation} 
Therefore, the quadratic form
$\bb_V[u]$ generates a bounded linear operator, say $\BB_V$, in the space
$\tH^1(\BG,x_0)$. It is easy to see that the operator $\BB_V$ is compact; 
actually, compactness will automatically follow from the estimates we
obtain in the next section. This operator is self-adjoint provided the
function $V$ is real-valued and it is non-negative provided $V\ge0$ a.e. As
usual, it is natural to identify the spectrum of the problem \eqref{1:1}
with the spectrum of the operator $\BB_V$. Still, we recall the
corresponding argument. 
\vskip0.2cm
 
The Laplacian $-\D$ on $\BG$, with the boundary conditions as in
\eqref{1:1},  is defined as the self-adjoint operator in
$\tL^2(\BG)$, associated with the quadratic form $\int_{\BG}|u'|^2dx$
considered on the domain $\tH^1(\BG,x_0)$. Given an element
$f\in\tL^2(\BG)$, the equality $-\D u=f$ 
under these boundary conditions means that $u$ is the unique
function in $\tH^1(\BG,x_0)$ such that
\begin{equation}\label{3:eq}
\int_{\BG}u'\overline{\varf'}dx=\int_{\BG}f\overline{\varf}dx,\qquad\forall
\varf\in\tH^1(\BG,x_0). 
\end{equation}
 
The Euler -- Lagrange equation reduces to $-u''=f$ on each edge. The 
continuity of $u$ on the whole of $\BG$ and the boundary
condition $u(x_0)=0$ are
implied by the inclusion
$u\in\tH^1(\BG,x_0)$. At each vertex $v\neq x_0$ the solution
$u$ meets the natural condition in the sense of Calculus of Variations.
Namely, let $e_1,\ldots, e_{d(v)}$ be the edges adjacent to a given vertex
$v$, oriented in such direction that $v$ is their initial point. Then
the condition at $v$, traditionally called {\sl Kirchhoff's condition}, is
\begin{equation*}
(u\res e_1)'(v)+\ldots+(u\res e_{d(v)})'(v)=0.
\end{equation*}
If, in particular, $d(v)=1$, this turns into $u'(v)=0$ which is the
boundary condition required by \eqref{1:1}.

The requirement $f\in\tL^2(\BG)$ is unnecessary for the existence of a
solution $u\in\tH^1(\BG,x_0)$ of the equation (or ``integral identity'')
\eqref{3:eq}. The solution
does exist if and only if $f$ is such that
the expression in the right-hand side generates
a continuous anti-linear functional on the space $\tH^1(\BG,x_0)$. One has
to take into account that the solution of \eqref{3:eq} for
$f\notin\tL^2(\BG)$ does not belong to the domain of the Laplacian
considered as the operator in $\tL^2(\BG)$. Still, it is conventional to
interpret \eqref{3:eq} as the weak form of the equation $-\D u=f$. In
particular, this is the case if $f\in\tL^1(\BG)$, due to the embedding
$\tH^1(\BG,x_0)\subset\tC(\overline{\BG})$. 
 
According to the above interpretation, the equation \eqref{1:1} means that
\begin{equation}\label{3:5}
\l\int_{\BG}u'\overline{\varf'}dx=\int_{\BG}Vu\overline{\varf}dx,\qquad\forall
\varf\in\tH^1(\BG,x_0). 
\end{equation} 
For $V\in\tL^1(\BG)$ we have also
$Vu\in\tL^1(\BG)$, so that the general scheme applies. 
\vskip0.2cm
 
On the other hand, the equation $\BB_Vu=\l u$ means that 
\begin{equation*}
\bb_V[u,\varf]:=\int_{\BG}Vu\overline \varf dx=
\l\int_{\BG}u'\overline{\varf'}dx,\qquad
\forall \varf\in\tH^1(\BG;x_0). 
\end{equation*} 
Comparing this with
\eqref{3:5}, we see that the eigenpairs for both equations are the same.
\vskip0.2cm
 
Below we use for the eigenvalues of \eqref{1:1} the notations
$\pm\l_n^{\pm}(\BB_V) $. If $V\ge 0$, we write $\l_n$ instead of $\l_n^+$. 
 
\subsection{Reduction to the case of trees} In order to have the
possibility to use \thmref{2:par}, it is necessary to reduce the problem
to the case of trees.  The procedure of such reduction (cutting cycles) is
rather standard. Nevertheless, we describe it in detail. 
 
\begin{lem}\label{3:red}
 Let $\BG$ be a compact metric graph. Then there
exist a compact metric tree $\BT$ and a continuous mapping $\tau:\BT\to\BG$
such that the operator $\tau^*:u(x)\mapsto u(\tau(x))$ defines an isometry
of the space $\tL^{1,2}(\BG)$ onto a
 subspace of finite codimension in $\tL^{1,2}(\BT)$.  
\end{lem} 
\begin{proof}
Let $\BG$ be a compact connected graph which is not a tree. Take any edge
$e_0=\lu v,w\ru\in\CE(\BG)$ which is part of a cycle in $\BG$,
then the subgraph $\BG_{e_0}=\BG\setminus\Int(e_0)$, with
$\CV(\BG_{e_0})=\CV(\BG)$ and $\CE(\BG_{e_0})=\CE(\BG)\setminus \{e_0\}$,
is connected. 
 
Choose a point $x\in\Int e_0$ and replace it by a pair $\{x_1,x_2\}$ of
new vertices. This gives rise to a new graph $\BG_1$ whose rigorous
definition is as follows.
Its sets of vertices and edges are defined by 

\begin{equation*} 
\CV(\BG_1)=\CV(\BG)\cup\{x_1,x_2\},\qquad
\CE(\BG_1)=\CE(\BG_{e_0})\cup \{\lu v,x_1\ru,\ \lu w,x_2\ru\}. 
\end{equation*} 
For any edge $e\in\CE(\BG_{e_0})$ its endpoints and its length
in $\BG_1$ are defined to be the same as in $\BG_{e_0}$. We set 
\begin{equation*} |\lu
v,x_1\ru|_{\BG_1}=|\lu v,x\ru|_{\BG},\qquad |\lu
w,x_2\ru|_{\BG_1}=|\lu w,x\ru|_{\BG}. 
\end{equation*} 
Evidently, $|\BG_1|=|\BG|$. 
 
Now we define a mapping $\tau_1:\BG_1\to\BG$.  Namely, $\tau_1$ is
identical on $\BG_{e_0}$ and
isometrically sends the edge $\lu v,x_1\ru$ onto $\lu v,x\ru$ and the edge
$\lu w,x_2\ru$ onto $\lu w,x\ru$.  
It is clear that the mapping $\tau_1$ is continuous and the
corresponding mapping $\tau_1^*:u(x)\mapsto u(\tau_1(x))$ is an isometry
of the space $\tL^{1,2}(\BG)$ into $\tL^{1,2}(\BG_1)$. 
Its image consists of all
functions $u\in\tL^{1,2}(\BG_1)$ such that $u(x_1)=u(x_2)$. Therefore, the
image is a subspace of codimension $1$. 
 
The number of simple cycles in $\BG_1$ is smaller than that for the
initial graph $\BG_0:=\BG$. Therefore, repeating this construction, we
obtain a finite sequence of graphs $\BG_0,\ldots, \BG_m$ such that the
graph $\BG_m$ has no cycles and loops (i.e. is a tree) and a sequence of
corresponding mappings $\tau_j:\BG_j\to \BG_{j-1}$, $j=1,\ldots, m$. The
tree $\BT=\BG_m$ and the mapping $\tau=\tau_1\circ\tau_2
\circ\ldots\circ\tau_m$ satisfy all the requirements of Lemma. 
\end{proof}
\vskip0.2cm It is useful to note also that 
\begin{equation*}
\int_{\BG}|u|^2Vdx=\int_{\BT}|\tau^*u|^2\tau^*Vdx 
\end{equation*} 
for any
$V\in\tL^1(\BG)$ and any $u\in\tC(\BG)$. 
 
 \section{Proof of the main results} First of all we prove \thmref{3:w}
and then derive \thmref{3:wid} from it. \thmref{1:eig} follows from the
latter almost immediately.

\subsection{ Proof of \thmref{3:w}.} 
1. Let first $\BG=\BT$ be a compact
tree, $\{T,\xi\}$ be its punctured subtree and $T=\T_1\cup\ldots\cup
\T_{d_T(\xi)}$ be the canonical partition of $\{T,\xi\}$. Applying the
inequality \eqref{3:m} to each subtree $\T_j$, we get for any
$u\in\tL^{1,2}(\BT)$: 
\begin{equation*} 
\int_{\T_j} |u(x)-u(\xi)|^2Vdx\le
|\T_j|\int_{\T_j}Vdx\,\int_{\T_j}|u'|^2dx, \qquad j=1,\ldots,d_T(\xi) 
\end{equation*} 
and hence, 
\begin{equation}\label{4:0}
\int_T|u(x)-u(\xi)|^2Vdx\le\int_T|u'|^2dx\max\limits_{j=1,\ldots,d_T(\xi)}
\biggl(|\T_j|\int_{\T_j}Vdx\biggr). 
\end{equation} 
Consider a function of
subtrees $T\subset\BT$, 
\begin{equation}\label{4:01}
\Phi(T)=\Phi(T;V)=|T|^{1/2}\bigl(\int_T Vdx\bigr)^{1/2}. 
\end{equation} 
By
\eqref{2:gen}, $\Phi\in\tS(\BT)$. The inequality \eqref{4:0} can be
written as 
\begin{equation}\label{4:03} 
\int_T |u(x)-u(\xi)|^2Vdx\le
\wt{\Phi}^2(T,\xi)\int_T|u'|^2dx. 
\end{equation} Here $\wt{\Phi}$ is the
function of punctured subtrees associated with the function \eqref{4:01},
cf.  \eqref{2:punc}.  \vskip0.2cm
 
Suppose now that the tree $\BT$ is split into the union of punctured
subtrees $\{T_1,x_1\},\ldots,\\ \{T_k,x_k\}$. Let $P$ be the corresponding
operator \eqref{3:op}. It follows from \eqref{4:03} that 
\begin{gather*}
\int_{\BT}|u-Pu|^2 Vdx=\sum_{j=1}^k \int_{T_j}|u(x)-u(x_j)|^2 Vdx\le
\sum_{j=1}^k \wt{\Phi}^2(T_j,x_j)\int_{T_j}|u'|^2dx  \\ \le
\|u'\|^2_2\,\bigl(\max\limits_j\wt{\Phi}(T_j,x_j)\bigr)^2,\qquad \forall
u\in\tL^{1,2}(\BT).\notag 
\end{gather*}
Now we apply \thmref{2:par}, to find that for any $n\in\N$ there exists a
partition of $\BT$ into $k$ punctured subtrees, such that $k\le n$ and
\begin{equation*} 
\int_{\BT}|u-Pu|^2 Vdx\le (n+1)^{-2}|\BT|\int_{\BT}
Vdx\,\|u'\|_2^2, \qquad\forall u\in\tL^{1,2}(\BT). 
\end{equation*} 
For the case of trees the proof is complete. 
 
2. Let $\BG$ be an arbitrary compact graph. Let $\BT$ and $\tau:\BT\to\BG$
be a compact tree and a mapping constructed according to Lemma \ref{3:red}. 
Without loss of generality, we assume that the pre-image $\tau^{-1}(x_0)$
consists of a single point. 
 
Let $\{T_j,x'_j\},\ j=1,\ldots,k$ be the partition of $\BT$, such that
\eqref{3:2} (with $\BT$ instead of $\BG$) is satisfied. Take
$G_j=\tau(T_j)$ and $x_j=\tau(x'_j)$. Since the mapping $\tau$ is continuous,
each set $G_j$ is closed and connected, and hence is a subgraph of $\BG$.
The punctured subgraphs
$\{G_j,x_j\}$ constitute a partition of $\BG$. Now an elementary 
computation shows that \eqref{3:2} is fulfilled for the graph $\BG$
if we take as $P$ the operator \eqref{3:op} corresponding to
this partition.

\subsection{Proof of \thmref{3:wid}.} 
By the
standard limiting argument, the inequality \eqref{3:m} extends from the
graph $\BG$ to its compactification $\overline{\BG}$. 
 
Let $\{\BG^{(m)}\}_{1\le m<\infty}$ be the family of compact subgraphs of
$\BG$, such that \eqref{2:su} is fulfilled. Fix a number $n\in\N$.
For any $m$, let 
\begin{equation*}
P_m:u\mapsto \sum_{j=1}^{k_m}
u(x^m_j)\chi^m_j
\end{equation*} 
be the
operators \eqref{3:op} for the subgraphs $\BG^{(m)}$ and the weight functions
$V\res \BG^{(m)}$. These operators depend also on $n$ but
we do not reflect this dependence  in the notations. 

According to the inequality \eqref{3:2}, we have for any $m$
and for any function $u\in\tL^{1,2}(\BG)$ normalized by the condition
$\int_{\BG}|u'|^2dx=1$: 
\begin{gather} \int_{\BG^{(m)}}|u-\sum_{j=1}^{k_m}
u(x^m_j)\chi^m_j|^2Vdx\label{4:inter}\\ \le
(n+1)^{-2}|\BG^{(m)}|\int_{\BG^{(m)}}Vdx\int_{\BG^{(m)}}|u'|^2dx \le
(n+1)^{-2}|\BG|\int_{\BG}Vdx.\notag 
\end{gather}
 
For different values of $m$, the numbers $k_m$ may be different. Denote by
$k=k(V)$ the minimal number such that $\#\{m:k_m=k\}=\infty$. Thinning
out the sequence $\{\BG^{(m)}\}$, we can assume that $k_m=k$ for all $m$.
Passing to a subsequence, we find points $\overline x_1,\ldots
\overline x_k\in\overline{\BG} $ (not necessarily different), such
that $x^m_j\to\overline x_j$
 as $m\to\infty$ for all $j=1,\ldots,k$. 
 
Now we show that the desired result is satisfied for the points found. 
Indeed, let $u\in\tL^{1,2}(\BG)$, $\int_{\BG}|u'|^2 dx=1$
 and $u(\overline x_1)=\ldots= u(\overline x_k)=0$. Then we have for each
$m$:  
\begin{equation}\label{4:3}
\biggl(\int_{\BG^{(m)}}|u|^2Vdx\biggr)^{1/2}\le \biggl(\int_{\BG^{(m)}}|u-P_m
u|^2Vdx\biggr)^{1/2} +\biggl(\int_{\BG^{(m)}}\bigl|\sum_{j=1}^k
u(x_j^m)\chi^m_j|^2 Vdx\biggr)^{1/2}. 
\end{equation} 
If $m$ is large
enough, then $|u(x_j^m)|=|u(x_j^m)-u(\overline x_j)|<\vare$ for all
$j=1,\ldots,k$ where $\vare>0$ is arbitrarily small. Then 
\begin{equation}\label{4:4} 
\int_{\BG^{(m)}}\bigl|\sum_{j=1}^k
u(x_j^m)\chi^m_j|^2Vdx \le \vare^2\int_{\BG^{(m)}}Vdx\le
\vare^2\int_{\BG}Vdx.  
\end{equation} 
Letting $m\to\infty$ in the
inequality \eqref{4:3} and taking \eqref{4:inter} and \eqref{4:4} into
account, we arrive at the desired inequality \eqref{3:wi}.\qed
 
\subsection{Proof of \thmref{1:eig}.} 
Let first $V\in\tL_+(\BG)$. Fix a
number $n\in\N$ and find points $\overline x_1,\ldots,\overline x_k$
according to \thmref{3:wid}. The subspace 
\begin{equation*}
\{u\in\tH^1(\BG,x_0): u(\overline x_1)=\ldots=u(\overline x_k)=0\}
\subset\tH^1(\BG,x_0)  
\end{equation*} 
is of codimension $k\le n$, and for
$n>1$ the inequality 
\begin{equation}\label{4:fin} 
\l_n(\BB_V)\le |\BG| n^{-2}\int_{\BG}Vdx,\qquad V\ge 0 
\end{equation} 
follows
from \eqref{3:wi} by the variational principle. The same inequality for
$n=1$ is implied by the estimate \eqref{3:r}. This completes the proof of
\eqref{1:ei} for $V\ge 0$. The result for sign-indefinite $V$ follows from
from \eqref{4:fin} by the variational principle, due to the inequalities
$\pm \bb_V[u]\le\bb_{V_{\pm}}[u]$. 
 \vskip0.2cm
The asymptotics \eqref{1:3} is an almost immediate consequence of the
estimate \eqref{1:ei}. Indeed, suppose first that the weight function $V$
has compact support, i.e. it vanishes outside a compact subgraph
$G\subset\BG$. Then we insert additional conditions $u(v)=0$ at all the
vertices $v\in \BG,\ v\neq x_0$. Since the number of these conditions is
finite, they do not affect the spectral asymptotics. For the new problem
the result follows from the well known asymptotic formula for a single
interval, and we are done. In the general case we fix $\vare>0$ and find a
compactly supported function $V_{\vare}$ such that
$\|V-V_{\vare}\|_1<\vare$. For the operator $\BB_{V_{\vare}}$ the
asymptotics \eqref{1:3} is satisfied, and for the operator
$\BB_V-\BB_{V_{\vare}}=\BB_{V-V_{\vare}}$ we have the estimate
\begin{equation*}
\l_n(|\BB_V-\BB_{V_{\vare}}|)\le|\BG|\vare n^{-2},\qquad\forall n\in \N. 
\end{equation*} 
Now the asymptotics \eqref{1:3} for the operator $\BB_V$
is implied by Lemma on continuity of the asymptotic coefficients, see
\cite{BS2}, Lemma 1.18.  \qed
 
\section{Complementary remarks. Applications to integral operators on
graphs} 
\subsection{On sharpness of the estimate \eqref{1:ei}.} 
Consider
the simplest case when $\BG$ is a single segment $[a,b]\subset\R$ and $V\ge0$.
The analog of \eqref{1:ei} for
the eigenvalue problem $-\L u''=Vu,\ u(a)=u(b)=0$ is the inequality 
\begin{equation}\label{5:dir} 
4n^2\L_n\le(b-a)\int_a^b Vdx. 
\end{equation} 
Its sharpness was discussed in
\cite{NS}, Section 3.2. It was shown there that for any $\vare>0$ and any
fixed $n_0\in\N$ a function $V=V_{\vare,n_0}$ can be found, such that the
corresponding eigenvalue $\L_{n_0}$ satisfies the inequality
\begin{equation*} 
4n_0^2\L_{n_0}\ge(1-\vare)(b-a)\int_a^b Vdx. 
\end{equation*} 
So, the constant in \eqref{5:dir} is sharp. 
   
Turn now to the eigenvalue problem $-\l u''=Vu,\ u'(0)=u(L)=0$, which is
just our problem \eqref{1:1} for the graph $\BG=[0,L]$, with $x_0=L$. 
Each eigenvalue $\l_n$ of this problem coincides with the eigenvalue
$\L_{2n-1}$ for the equation $-\l u''=V(|x|)u$ on the interval $(-L,L)$,
with the zero boundary conditions at both ends. Using the above result, we
find a function $V\ge0$ such that 
\begin{gather*}
(2n_0-1)^2\l_{n_0}=(2n_0-1)^2\L_{2n_0-1}\ge (1-\vare)L\int_0^L Vdx. 
\end{gather*} 
Taking here $n_0=1$, we see that constant factor $1$ in
\eqref{1:ei} is the best possible. However, for each particular $n$ the
factor $1$ is probably not sharp. 
\vskip0.2cm 
Note also that the
inequality \eqref{3:m} implies the estimate 
\begin{equation*} 
\|\BB_V\|\le
\diam(\BG)\int_{\BG}|V|dx,\qquad \diam(\BG):=\sup\{\rho(x,y):x,y\in\BG\}.
\end{equation*} 
It follows that the inequality \eqref{1:ei} can be
replaced by 
\begin{equation*}
\l_n^{\pm}\le\min\biggl(\frac{|\BG|}{n^2},\diam(\BG)\biggr)
\int_{\BG}V_{\pm}dx. 
\end{equation*} 
This can be useful when dealing with graphs of
small diameter but large total length. 
 
\subsection{Singular numbers of the operator $a(x)(-\D)^{-1/2}$ in
$\tL^2(\BG)$.} Recall that the singular numbers ($s$-numbers) of a compact
operator $B$ acting between two Hilbert spaces, are defined as the
non-zero eigenvalues of any of two self-adjoint compact operators
$|B|:=(B^*B)^{1/2}$ and $|B^*|:=(BB^*)^{1/2}$. \vskip0.2cm
 
According to the
Hilbert space theory, for each $u\in\tH^1(\BG,x_0)$ we have
$\|u'\|_2=\|(-\D)^{1/2}u\|_2$. Here $-\D$ is the operator whose rigorous
description was given in Subsection 3.2.
For any $V\in\tL_+(\BG)$ and for
any non-zero element $u\in\tH^1(\BG,x_0)$ the equality holds
\begin{equation}\label{5:rat}
\frac{\int_{\BG}|u|^2Vdx}{\int_{\BG}|u'|^2dx}=
\biggl(\frac{\|V^{1/2}(-\D)^{-1/2}w\|_2}{\|w\|_2}\biggr)^2, \qquad
w=(-\D)^{1/2}u. 
\end{equation} 
Taking into account the variational description of
the eigenvalues and the $s$-numbers, we conclude from \eqref{5:rat} that
\begin{equation*}
\l_n(\BB_V)=s_n^2\bigl(V^{1/2}(-\D)^{-1/2}\bigr),\qquad\forall n\in\N. 
\end{equation*} 
Let now $a(x)$ be an arbitrary function from $\tL^2(\BG)$ and
$V(x)=|a(x)|^2$, then $a(x)=\psi(x) V^{1/2}(x)$ where $|\psi(x)|=1$ a.e. 
Multiplication by 
$\psi$ is a unitary operator.  Therefore, the
operators $V^{1/2}(x)(-\D)^{-1/2}$ and $a(x)(-\D)^{-1/2}$ are compact
simultaneously, and have the same $s$-numbers. This leads to a useful
consequence of the estimate \eqref{1:ei}, 
 
\begin{equation}\label{5:x} 
s_n\bigl(a(x)(-\D)^{-1/2}\bigr)\le
\frac{|\BG|^{1/2}\|a\|_2}{n}, \ \forall n\in\N;\qquad a\in\tL^2(\BG).
\end{equation} 
 
\subsection{The case $V\equiv 1$.} 
In this case $\BB_V=(-\D)^{-1}$, thus $\l_n(\BB_V)=\l_n^{-1}(-\D)$ and the
relations \eqref{1:ei} and \eqref{1:3} turn into 
\begin{gather}
|\BG|^2\l_n(-\D)\ge n^2,\qquad\forall n\in\N;\label{5:ei}\\
\frac{\sqrt{\l_n(-\D)}}{n}\to\frac{\pi}{|\BG|}.\label{5:3} 
\end{gather}
Even for this simplest case, the estimate \eqref{5:ei} is informative,
since it is uniform with respect to all graphs of a given length,
independently of their combinatorial structure. By means of standard
perturbation arguments, the asymptotics \eqref{5:3} extends to the
operators $-\D+q$ of Sturm -- Liouville type with the real-valued, bounded
potential $q$. This considerably improves the result of Theorem 5.4 of the
paper \cite{C}. For such operators on the so-called regular trees, another
proof of the asymptotics \eqref{5:3} was recently given in \cite{S}. 
 
\subsection{Estimates of singular numbers for integral operators on
graphs} The above results admit immediate applications to
the estimates of singular numbers of the integral operators on graphs.
We restrict ourselves with the simplest case when the operator acts in 
$\tL^2(\BG)$ according to the formula
\begin{equation}\label{5:int}
\bigl(\BK u\bigr)(x)=\int_{\BG}K(x,y)u(y)dy.  
\end{equation}
We suppose that the kernel $K(x,y)$ belongs to the space
$\tL^2_y\bigl(\BG,\tW^{1,2}_x(\BG)\bigr)$. This means that 
$K$ is measurable on $\BG\times\BG$, for  
almost all $y\in\BG$ the function $K(\cdot,y)$ lies in the space
$\tL^{1,2}(\BG)$ and, moreover,
\begin{equation*}
\cm(K,\BG):=\int_{\BG}\int_{\BG}\bigl(|K(x,y)|^2+
|\BG|^2|K'_x(x,y)|^2\bigr)dxdy <\infty. 
\end{equation*} 
The expression for
$\cm(K,\BG)$ is homogeneous with respect to the similitudes of the graph
$\BG$. This explains why the factor $|\BG|^2$ is included.

\begin{thm}\label{5:sinn} Let $\BG$ be a graph of finite total
length and 
$K\in\tL^2_y\bigl(\BG,\tW^{1,2}_x(\BG)\bigr)$. Then the following
inequality is satisfied for the $s$-numbers of the operator
$\BK$ given by \eqref{5:int}: 
 
\begin{equation}\label{5:z} 
\sum_{n=1}^{\infty}n^2 s^2_n(\BK)\le
32|\BG|^2 \cm(K,\BG). 
\end{equation}
 
If in addition there exists a point $x_0\in\BG$ such that $K(x_0,y)=0$ for
almost all $y\in\BG$, then the inequality \eqref{5:z} can be refined:
\begin{equation}\label{5:zz} 
\sum_{n=1}^{\infty}n^2 s^2_n(\BK)\le
8|\BG|^2 \int_{\BG}\int_{\BG}|K'_x(x,y)|^2dxdy. 
\end{equation} 
\end{thm} 
\begin{proof} 
We follow the approach of \cite{GK}, Section II.10.4.
\vskip0.2cm

We start with the proof of the second statement of Theorem. Its assumption
means that the function $K(\cdot,y)$ lies
in $H^1(\BG,x_0)$ for almost all $y\in\BG$. This allows one to write
\begin{equation*} 
K(\cdot,y)=(-\D_x)^{-1/2}L(\cdot,y)
\end{equation*}
where
$L(\cdot,y)=(-\D_x)^{1/2}K(\cdot,y)$ and 
$-\D_x$ is the
operator $-\D$ (with the boundary conditions as in \eqref{1:1}, cf.
Subsection 3.2), acting in the variable $x$. This representation
of the kernel yields factorization of the corresponding
operator, $\BK=(-\D_x)^{-1/2}\BL$. 

From the estimate \eqref{5:x} for $a\equiv 1$ we derive 
$s_n\bigl((-\D_x)^{-1/2}\bigr)\le |\BG|n^{-1}$.
The operator $\BL$ is of the Hilbert-Schmidt class and moreover,
\begin{equation}\label{5:6} 
\sum_{n=1}^{\infty}s_n^2(\BL)=\int_{\BG}\int_{\BG}|L(x,y)|^2 dxdy=
\int_{\BG}\int_{\BG} |K'_x(x,y)|^2 dxdy. 
\end{equation}
Due to an inequality by Ky Fan, see \cite{GK},
Corollary II.2.2, 
\begin{equation*} 
s_{2n}(\BK)\le
s_{2n-1}(\BK)\le s_n\bigl((-\D)^{-1/2}\bigr) s_n(\BL)\le
|\BG|n^{-1}s_n(\BL). 
\end{equation*} 
Together with \eqref{5:6} this implies \eqref{5:zz}.

\vskip0.2cm
Let us turn to the proof of \eqref{5:z}.
The function $K(\cdot,y)$ is
continuous on $\BG$ for almost all $y\in\BG$. Choose a point $x_0\in\BG$
and split the kernel $K$ into the sum $K=K^0+\wt K$ where
$K^0(x,y)=K(x_0,y)$. The operator $\BK^0$ has rank one, and its only
non-zero singular number is given by 
\begin{equation}\label{5:zy}
s_1^2(\BK^0)=|\BG|\int_{\BG}|K(x_0,y)|^2dy \le
2\cm(K,\BG). 
\end{equation}
Let us explain the latter inequality. Any function $w\in\tL^{1,2}(\BG)$ 
assumes its mean value $\hat w=|\BG|^{-1}\int_{\BG}w dx$ at some point
$\xi\in\overline{\BG}$. Using the inequality \eqref{3:m} (which extends
from $\BG$ to any $\xi\in\overline{\BG}$ by the continuity), we obtain
$|w(x_0)-\hat w|^2\le |\BG|\int_{\BG}|w'|^2 dx$.
Also,
$|\BG||\hat w|^2\le \int_{\BG}|w|^2 dx$.
Therefore, 
$|\BG||w(x_0)|^2\le 2\int_{\BG}\bigl(|w|^2+|\BG|^2 |w'|^2\bigr)dx$,
whence \eqref{5:zy}.

Since $\BK^0$ is a rank one operator, one has 
$s_{n+1}(\BK)\le s_n(\wt{\BK}),\ \forall n\in \N$.
For $\wt{\BK}$ the inequality \eqref{5:zz} is satisfied, because 
$\wt K(x_0,y)=0$
a.e.
Therefore, 
\begin{equation}\label{5:5}
\sum_{n=2}^{\infty}n^2s^2_n(\BK)\le\sum_{n=1}^{\infty}(n+1)^2s^2_n(\wt\BK)
\le 32|\BG|^2
\int_{\BG}\int_{\BG}|K'_x(x,y)|^2 dxdy. 
\end{equation} 
Further,
\eqref{5:zz} implies $s_1(\wt{\BK})\le\sqrt 8 \cm(K,\BG)$. 
The inequality \eqref{5:z} follows from here,
\eqref{5:5} and
\eqref{5:zy}, due to the triangle inequality $s_1(\BK)\le
s_1(\wt{\BK})+s_1(\BK^0)\le 3\sqrt2\,\cm(K,\BG)$. 
\end{proof}
\begin{cor}\label{5:indiv}  
Under the assumptions of \thmref{5:sinn} the
singular numbers $s_n(\BK)$ of the operator \eqref{5:int} satisfy
the estimate 
\begin{equation}\label{5:fin}
s_n(\BK)\le\frac{4\sqrt 6}{n^{3/2}}\cm^{1/2}(K,\BG), \qquad \forall
n\in\N 
\end{equation} 
and besides, $s_n(\BK)=o(n^{-3/2})$. 
\end{cor} 
\begin{proof} 
Denote $s_n=s_n(\BK)$ and $C^2=
32\cm(K,\BG)$. The estimate \eqref{5:fin} is implied by
the inequality 
\begin{equation*} \frac{n^3}{3}s_n^2\le s_n^2\sum_{k=1}^n
k^2 \le \sum_{k=1}^n k^2 s_k^2\le C^2,\qquad \forall n\in\N. 
\end{equation*} 
The relation $s_n(\BK)=o(n^{-3/2})$ follows from the
inequality 
\begin{equation*} 
c n^3 s_n^2\le \sum_{k=[n/2]}^n k^2s_k^2,
\qquad n>1,\ c>0 
\end{equation*} 
in which the right-hand side tends to
zero as $n\to\infty$, due to the convergence of the series in \eqref{5:z}.
\end{proof}

We would like to emphasize that the estimates \eqref{5:z}, \eqref{5:zz}
and \eqref{5:fin}
are uniform with respect to all graphs of a given length. Specific values
of the constants in these estimates are not so important.

In the same way, it is possible to study similar operators (with $dy$ in
\eqref{5:int} replaced by $d\mu(y)$) acting between the spaces $\tL^2(\BG,\mu)$
and $\tL^2(\BG,\nu)$, where $\mu$ and $\nu$ are finite Borelian measures.
Such operators appear in various applications, see the paper \cite{BS3}.
Note also that for the case of trees an analog of \thmref{5:sinn} for $l$ 
times differentiable
kernels can be derived from \thmref{6:eigh} of the next section. 
 
\section{Spaces $\tH^l(\BT,x_0)$ and operators $\BB_{l,V}$.}
 
Let us discuss the higher order analogs of the space $\tL^{1,2}$.  
Here a
serious obstacle arises, since the continuity of derivatives
$u',\ldots,u^{(l-1)}$ at the vertices should be included in the definition. 
However, the derivatives of odd order change their sign depending on the
orientation on edges, so that for $l>1$ the space $\tL^{l,2}(\BG)$ can be well
defined only for oriented graphs, and for different choices of orientation
such spaces are substantially different. For this reason, we define the
spaces $\tL^{l,2}$ only on trees, since for them a natural orientation 
does exist. 
 
So, let $\BG=\BT$ be a tree of finite total length
and let $x_0\in\BT$ be a vertex selected (the
root).  The natural
partial ordering on the rooted tree $\{\BT,x_0\}$ is introduced as
follows: 
\begin{equation*} 
x\preceq y\Longleftrightarrow x\in\lu x_0,
y\ru. 
\end{equation*} 
Recall that $\lu x_0, y\ru$ is the unique simple
polygonal path in $\BT$ connecting $x_0$ with $y$.
We always parametrize the edges of $\BT$ in the
direction, compatible with this partial ordering. 
 
Now we are in a position to define the space $\tH^l(\BT,x_0)$ for
arbitrary $l\in\N$. A function $u$ on $\BT$ belongs to $\tH^l(\BT,x_0)$ if
$u$ is continuous on $\BT$, the restriction of $u$ to each edge $e$ lies
in $\tH^l(e)$, the functions $u',\ldots u^{(l-1)}$ extend from
$\BT\setminus\CV(\BT)$ to the whole of $\BT$ as continuous functions,
$u(x_0)=\ldots=u^{(l-1)}(x_0)=0$ and $u^{(l)}\in \tL_2(\BT)$.
We consider $\tH^l(\BT,x_0)$ as the
Hilbert space with the scalar product $(u,v)_{\tH^l(\BT,x_0)}=
(u^{(l)},v^{(l)})_{\tL^2(\BT)}$ and the corresponding norm. 
\vskip0.2cm

Let $\xi,x\in\BT$ and $\xi\preceq x$.  Consider the Taylor polynomial
\begin{equation*}
P_{l-1}(t;u,\xi)=\sum_{k=0}^{l-1}\frac{u^{(k)}(\xi)t^k}{k!},
\end{equation*} 
then, due to our agreement about orientation,
\begin{equation*}
u(x)-P_{l-1}(\rho(\xi,x);u,\xi)=\frac{1}{(l-1)!}\int_{\lu\xi,x\ru}u^{(l)}(y) 
\rho^{l-1}(y,x)dy. 
\end{equation*} 
By Cauchy's inequality,
\begin{equation}\label{6:m} 
((l-1)!)^2 |u(x)-P_{l-1}(\rho(\xi,x);u,\xi)|^2
\le\frac{\rho^{2l-1}(\xi,x)}{2l-1}\int_{\lu\xi,x\ru}|u^{(l)}(y)|^2dy.
\end{equation}
 
Given a function $V\in\tL^1(\BT)$, let $\bb_V$ is the corresponding
quadratic form, cf. \eqref{3:s}. It follows from the inequality
\eqref{6:m} for $\xi=x_0$ that 
\begin{gather*} 
|\bb_V[u]|\le
C'(l)|\BT|^{2l-1}\|u^{(l)}\|^2_2\int_{\BT}|V(x)|dx,\qquad\forall
u\in\tH^l(\BT,x_0), \\ C'(l)=((l-1)!)^{-2}(2l-1)^{-1}.
\end{gather*} 
Therefore, the quadratic form $\bb_V[u]$ generates in
$\tH^l(\BT,x_0)$ a bounded linear operator which we denote by $\BB_{l,V}$.
Its eigenpairs $\{\l,u\}$ correspond to the problem 
\begin{gather} \l(-\D)^lu=Vu,\qquad u(x_0)=u'(x_0)=\ldots=u^{(l-1)}(x_0)=0;
\label{6:1h}\\
u^{(l)}(v)=\ldots=u^{(2l-1)}(v)=0\ \text{if}\ v\in\p\BT\setminus\{x_0\}.
\notag 
\end{gather}

For rooted trees, \thmref{1:eig} is a particular case of the following
result.
 
\begin{thm}\label{6:eigh} 
Let $\BT$ be a rooted tree of finite total
length, $x_0$ be its root, and let $V=\overline
V\in\tL^1(\BT)$. Then the eigenvalues of the problem \eqref{6:1h} satisfy
the inequality 
\begin{equation}\label{6:eih} 
\l_n^{\pm}\le
C(l)\frac{|\BT|^{2l-1}\int_{\BT} V_{\pm}dx}{n^{2l}},\qquad
C(l)=\frac{l^{2l}}{(l-1)!)^2(2l-1)},\qquad \forall n\in\N. 
\end{equation}
Along with the estimate \eqref{6:eih}, the Weyl-type asymptotics holds,
\begin{equation*}
n\,\bigl(\l_n^{\pm}\bigr)^{\frac{1}{2l}}\to\pi^{-1}\int_{\BT}
\bigl(V_{\pm}(x)\bigr)^{\frac{1}{2l}}dx,\qquad n\to\infty. 
\end{equation*}
\end{thm} 
\noindent{\it Outline of the proof.} We discuss only 
the estimate \eqref{6:eih} for the case of
compact trees and $V\in\tL_+(\BT)$, since the rest needs no serious changes
compared with the proof of Theorems \ref{3:wid} and \ref{1:eig}. We also
suppose that $d(x_0)=1$. Otherwise, the operator $\BB_{l,V}$ splits
into the orthogonal sum of similar operators for each of $d(x_0)$ 
subtrees constituting the canonical partition of the punctured
tree $\{\BT,x_0\}$, cf. Section 2. The estimate \eqref{6:eih} for $\BB_{l,V}$
easily follows from the same estimate for its corresponding parts.
 
Consider a function of subtrees $T\subset\BT$: 
\begin{equation}\label{6:0} 
\Phi(T)=|T|^{1-\frac{1}{2l}}\biggl(\int_T
Vdx\biggr)^{\frac{1}{2l}}, 
\end{equation} 
and let $\wt\Phi(T,x)$ be the
function of punctured subtrees associated with $\Phi$, cf. \eqref{2:punc}. 
Suppose
that the punctured subtrees $\{T_j,x_j\}, j=1, \ldots,k$ constitute a
partition of $\BT$ and let $\chi_j$ be the characteristic function of
$T_j$. It can be easily derived from \eqref{6:m} that for
$u\in\tH^{l}(\BT,x_0)$ the inequality 
\begin{equation*}
\int_{\BT}\bigl|u-\sum_{j=1}^k
P_{l-1}(\rho(x_j,x);u,x_j)\chi_j\bigr|^2Vdx \le
C'(l)\|u^{(l)}\|^2_2\,\bigl(\max\limits_{j=1,\ldots,k}
(\wt\Phi(T_j,x_j)\bigr)^{2l}, 
\end{equation*} 
holds, provided all the
subtrees $\{T_j,x_j\}$ are oriented coherently to the orientation on
$\BT$, that is if $x_j\preceq x$ for any $x\in T_j$. 

\thmref{2:par} applies to the function $\Phi$ but this does not lead
automatically to the inequality \eqref{6:eih}. Indeed, we have to check
that all the subtrees $\{T_j,x_j\}$ are properly oriented. For this
purpose, let us return to Lemma \ref{1:lem} whose consequence is
\thmref{2:par}. The proof of Lemma started with choosing a vertex 
$v_0\in\p\BT$, then a path $\CP$ and subtrees $T_x^+$ for each $x\in\CP$
were constructed. The assumption
$d(x_0)=1$ means that $x_0\in\p\BT$.  If we take
$v_0=x_0$, then all the subtrees $T_x^+$ are oriented coherently to the 
orientation on $\BT$,
and the scheme goes through.

Applying inequality \eqref{1:main} to the function $\Phi$ introduced by
\eqref{6:0} and using the variational principle, we find that
\begin{equation*} \
\l_{Nl+1}(\BB_{l,V})\le
C'(l)\frac{|\BT|^{2l-1}\int_{\BT}Vdx}{N^{2l}}, \ \forall N\in\N;\qquad
V\in L_+(\BT). 
\end{equation*} 
The estimate \eqref{6:eih} for all $n\in\N$,
with $C(l)=C'(l)l^{2l}$, follows from here due to the monotonicity in $n$
of eigenvalues. 
 
\bibliographystyle{amsplain}

\end{document}